\documentclass[12pt]{article}
\begin{document}

\def\2{\frac{1}{2}}
\def\om{\omega}
\def\omt{\tilde{\omega}}
\def\ti{\tilde}
\def\o{\Omega}
\def\bchi{\bar\chi^i}
\def\In{{\rm Int}}
\def\ba{\bar a}
\def\w{\wedge}
\def\ep{\epsilon}
\def\k{\kappa}
\def\Tr{{\rm Tr}}
\def\ST{{\rm STr}}
\def\ss{\subset}
\def\ot{\otimes}
\def\bc{{\bf C}}
\def\br{{\bf R}}
\def\de{\delta}
\def\tr{\triangleleft}
\def\al{\alpha}
\def\la{\langle}
\def\ra{\rangle}
\def\G{\Gamma}
\def\th{\theta}
\def\lm{\ti\lambda}
\def\U{\Upsilon}
\def\jp{{1\over 2}}
\def\js{{1\over 4}}
\def\d{\partial}

\def\be{\begin{equation}}
\def\ee{\end{equation}}
\def\bea{\begin{eqnarray}}
\def\eea{\end{eqnarray}}
\def\D{{\cal D}}
\def\G{{\cal G}}
\def\H{{\cal H}}
\def\R{{\cal R}}
\def\T{{\cal T}}
\def\bT{\bar{\cal T}}
\def\F{{\cal F}}
\def\n{{1\over n}}
\def\si{\sigma}
\def\ta{\tau}
\def\ov{\over}
\def\l{\lambda}
\def\L{\Lambda}

\def\pih{\hat{\pi}}
\def\Vt{V^{\ti}}
\def\Ut{U^{\ti}}

\def\e{\varepsilon}
\def\b{\beta}
\def\ga{\gamma}

\begin{titlepage}
\begin{flushright}
{}~
FSU-TPI 03/01\\
\end{flushright}

\vspace{3cm}
\begin{center}
{\Large \bf WZW--Poisson manifolds}\\
[50pt]{\small
{\bf C. Klim\v{c}\'{\i}k}
\\ ~~\\Institute de math\'ematiques de Luminy,
 \\163, Avenue de Luminy, 13288 Marseille, France
\vskip1pc
and
\vskip1pc
{\bf T. Strobl}
\\ Institut f\"ur  Theoretische Physik \\Friedrich Schiller Universit\"at,
 \\ Max--Wien--Platz 1, D--07743 Jena, Germany}

\vspace{1cm}
\begin{abstract}
 We observe that a term of the WZW-type can be added to the Lagrangian of
 the Poisson $\si$-model in such a  way that the algebra of the first class
 constraints remains closed. This leads to a natural generalization of the
 concept of Poisson geometry. The resulting "WZW--Poisson" manifold $M$ is
 characterized by a bivector $\Pi$ and by a closed three-form $H$ such that
 $\2 [\Pi,\Pi]_{Schouten}\; = \; \langle H,\Pi\otimes\Pi\otimes\Pi \rangle$.
\end{abstract}
\end{center}
 \end{titlepage}
\noindent {\bf 1.} It turns out to be a  fruitful idea to associate dynamical
systems---and in particular topological field theories---to geometric data on
manifolds. Here we shall study the following example: Given a bivector
$\Pi=\jp\Pi^{ij}\partial_i\w\partial_j$ and a two-form
$\Omega=\jp\Omega_{ij}dX^i\w dX^j$
 on a manifold $M$, we can immediately write
down the action functional \be S[X,A]= \int_\Sigma (A_i \wedge dX^i + \jp
\Pi^{ij}(X) A_i \wedge A_j+\jp\Omega_{ij}(X)dX^i\w dX^j). \ee In the story
that follows, $\Sigma$ will be a cylindrical world-sheet, $X^i$ is a 
collection of coordinates on the target space $M$, and $A_i$ is a set of
$1$-forms on $\Sigma$. Of course this action can be written also in a
coordinate-independent way.

Introducing the standard world-sheet coordinates $\si$ and $\tau$ (the loop
and the evolution parameters, respectively) we set \be A_i=A_{i\si}d\si
+A_{i\tau}d\tau\ee and rewrite (1) in the following form \be
S[X,A]=\int_\Sigma d\si d\tau[(p_{i} \d_\tau{X^i}-A_{i\tau}\phi^{i}].\ee
Here \be p_i=A_{i\si}-\Omega_{ij}\d_\si X^j, \quad \phi^{i}=\d_\si
X^i+\Pi^{ij}p_{j} +\Pi^{ij}\Omega_{jk}\d_\si X^k.\ee

\noindent {\bf 2.} Let $P$ be some (possibly infinite-dimensional) manifold
equipped with a symplectic form $d\theta$. Suppose there is a set of
functions $h,\phi^\al,d^{\al}_{\b},c^{\al\b}_{\ga}$ fulfilling \be
\{h,\phi^\al\}=d^{\al}_{\b}\phi^\b,\quad
\{\phi^\al,\phi^\b\}=c^{\al\b}_{\gamma}\phi^\ga,\ee where the indices take
values in some set $U$ and the Poisson bracket corresponds to $d\theta$. To
these data we associate the constrained dynamical system described by the
action \be S= \int (\theta -(h +\lambda_\al\phi^\al)d\tau),\ee where
$\lambda_\al$ is a set of Lagrange multipliers and $\phi^\al$ are
the corresponding first class constraints.

Now the question arises: for which pair $\Pi,\Omega$ the action (3) defines
a constrained dynamical system in the sense described above (i.e. the
relations (5) should hold). Of course, $h=0$, $A_{i\tau}$ play the role of
the Lagrange multipliers $\lambda_\al$ and $\theta=\oint p_i d X^i$.

It is simple to  answer  this question. The symplectic form $d\theta$ has
the canonical Darboux
form and the calculation of the Poisson brackets is straightforward. We obtain
\be \{\phi^i(\si),\phi^j(\si')\}= -(\d_k\Pi^{ij}+
\Pi^{il}\Pi^{jm}(d\Omega)_{klm})\delta(\si-\si')\phi^k(\si),\ee
provided
\be \2 [\Pi,\Pi]_{S} \; = \;
\langle d\Omega,\Pi\otimes\Pi\otimes\Pi \rangle  \label{condition} \ee
holds true. The symbol $[.,.]_S$ denotes the Schouten bracket and the functions
$c^{\al\b}_{\ga}$ can be read off from (7). The contraction on the right hand
side is with respect to the first, third and fifth entry of $\Pi^3$. We remark
that the condition (\ref{condition}) is necessary and sufficient for the
system of constraints following from (1) to be of the first class (cf.\
\cite{inprep} for further details).

\noindent {\bf 3.} Our discussion can be slightly generalized. Consider the
bivector $\Pi$ and a {\it closed} $3$-form on the manifold $M$. To these
data we associate the following action \be S[X,A]= \int_\Sigma (A_i \wedge
dX^i + \jp \Pi^{ij}(X) A_i \wedge A_j) +\int_V H. \ee Here $V$ is the
interior   of the cylinder $\Sigma$ and by $H$ we really mean the pullback of
$H$ to $V$ by an extension to $V$ of the map $X^i(\si,\tau)$.  Of course,
there are the subtleties concerning the boundaries of the cylinder and the
WZW term. We do not give the detailed discussion in this letter. It is
a straightforward generalization of the treatment in \cite{KS},   where the
WZW model  on the cylinder is studied from  the point of view of
Hamiltonian mechanics.

Note that (9) reduces to (1) for $H=d\Omega$.
By repeating the previous discussion, we arrive at the conclusion that the
model (9) corresponds to a maximally constrained dynamical system iff \be
\2 [\Pi,\Pi]_{S}\; =\; \langle H,\Pi\otimes\Pi\otimes\Pi \rangle\;.\ee

\noindent {\bf 4.} For $H=0$, the action (9) defines the Poisson
$\si$-model \cite{SS,I,ASS} and the condition (10) says that the bivector
$\Pi$ satisfies the Jacobi identity. Therefore Poisson geometry could
have been invented by asking the question when the model (9) (with $H=0$)
is a maximally constrained dynamical system or a topological field theory.
If we do not set $H=0$, the same logic gives a natural generalization: the
concept of what one might call WZW-Poisson manifolds.
We repeat that the latter is characterized by a bivector $\Pi$ and a closed
$3$-form $H$ such that the condition (10) holds.

It remains to understand the properties of the WZW-Poisson manifolds in
more detail. It may be
that there is a non-trivial intersection of this notion with the other
generalizations of Poisson geometry like quasi-Poisson manifolds
\cite{AKM}, Dirac manifolds \cite{C}
or the manifolds leading to the nonassociative generalization \cite{CS,HM}
of the Kontsevich
expansion.

\noindent {\bf Note added:} After completion of this work we became
aware that the relation (10) was obtained also in \cite{Park} within a
BV approach.

  \end{document}